\documentclass[12pt]{article}
\pdfoutput=1

\usepackage{amsmath,amsfonts,amssymb,amsthm,bbm,color,graphicx,url}
\usepackage{natbib}
\usepackage[a4paper]{geometry}
\usepackage{hyperref}

\theoremstyle{plain}
\newtheorem{theorem}{Theorem}[]

\theoremstyle{definition}
\newtheorem{example}[]{Example}
\newtheorem{definition}[]{Definition}

\theoremstyle{remark}

\newcommand{\dd}{\,{\rm d}}
\newcommand{\one}[1]{\ensuremath{\mathbbm{1}}(#1)}

\newcommand{\real}{\mathbb{R}}
\newcommand{\Real}{\bar{\mathbb{R}}}

\newcommand{\cA}{\mathcal{A}}
\newcommand{\cB}{\mathcal{B}}
\newcommand{\cL}{\mathcal{L}}
\newcommand{\cP}{\mathcal{P}}
\newcommand{\cX}{\mathcal{X}}

\newcommand{\myL}{{\rm L}}
\newcommand{\myS}{{\rm S}}

\newcommand{\crps}{\mathrm{CRPS}}

\allowdisplaybreaks

\begin{document}

\begin{center}

\bf \Large 

Properization:

Constructing Proper Scoring Rules via Bayes Acts

\end{center}

\centerline{\bf Jonas R. Brehmer}
\centerline{University of Mannheim, Mannheim, Germany}

\bigskip

\centerline{\bf Tilmann Gneiting}
\centerline{Heidelberg Institute for Theoretical Studies, Heidelberg, Germany} 
\centerline{Karlsruhe Institute of Technology, Karlsruhe, Germany} 

\bigskip

\centerline{August 16, 2018}

\medskip

\begin{abstract} 

Scoring rules serve to quantify predictive performance.  A scoring
rule is proper if truth telling is an optimal strategy in expectation.
Subject to customary regularity conditions, every scoring rule can be
made proper, by applying a special case of the Bayes act construction
studied by \citet{GruenDawid2004} and \citet{Dawid2007}, to which we
refer as properization.  We discuss examples from the recent
literature and apply the construction to create new types, and
reinterpret existing forms, of proper scoring rules and consistent
scoring functions. In an abstract setting, we formulate sufficient
conditions under which Bayes acts exist and scoring rules can be made
proper.
 

\end{abstract}

\section{Introduction}  \label{sec:introduction}

Let $\cB$ be a $\sigma$-algebra of subsets of a general sample space
$\Omega$.  Let $\cP$ be a convex class of probability measures on
$(\Omega, \cB)$.  A {\em scoring rule}\/ is any extended real-valued
function $\myS$ on $\cP \times \Omega$ such that 
\[
\myS(P,Q) = \int \myS(P,\omega) \dd Q(\omega) 
\]
is well-defined for $P, Q \in \cP$.  The scoring rule $\myS$ is {\em
  proper}\/ relative to $\cP$ if
\begin{equation}  \label{eq:proper} 
\myS(Q,Q) \leq \myS(P,Q) \quad \mbox{for all} \quad P, Q \in \cP.   
\end{equation}
In words, we take scoring rules to be negatively oriented penalties
that a forecaster wishes to minimize.  If she believes that a future
quantity or event has distribution $Q$, and the penalty for quoting
the predictive distribution $P$ when $\omega$ realizes is
$\myS(P,\omega)$, then~\eqref{eq:proper} implies that quoting $P = Q$
is an optimal strategy in expectation.  The scoring rule is {\em
  strictly proper} if~\eqref{eq:proper} holds with equality only if $P
= Q$.  For recent reviews of the theory and application of proper
scoring rules see \citet{Dawid2007}, \citet{GneitRaft2007}, 
\citet{DawidMusio2014}, and \citet{GneitKatz2014}.

The intent of this note is to draw attention to the simple fact that,
subject to customary regularity conditions, any scoring rule can be
{\em properized}, in the sense that it can be modified in a
straightforward way to yield a proper scoring rule, so that truth
telling becomes an optimal strategy.  Implicitly, this construction
has recently been used by various authors in various types of
applications; see, e.g., \citet{Diks2011}, \citet{Christetal2014} and
\citet{HolzKlar2018}.

\begin{theorem}[properization]  \label{th:properization} 
Let\/ $\myS$ be a scoring rule.  Suppose that for every\/ $P \in \cP$
there is a probability distribution\/ $P^* \in \cP$ such that 
\begin{equation}  \label{eq:Bayes} 
\myS(P^*,P) \leq \myS(Q,P) \quad \mbox{for all} \quad Q \in \cP. 
\end{equation} 
Then the function
\begin{equation}  \label{eq:S*} 
\myS^* : \cP \times \Omega \to \Real, \quad 
(P,\omega) \mapsto \myS^*(P,\omega) = \myS(P^*,\omega),  
\end{equation} 
is a proper scoring rule. 
\end{theorem} 

Here and in what follows we denote the real line by $\real$ and the
extended real line by $\Real := \real \cup \{ - \infty, \infty \}$.
Any probability measure $P^*$ with the property~\eqref{eq:Bayes} is
commonly called {\em Bayes act}; for the existence of Bayes acts,
see Section~\ref{sec:existence}.  In case there are multiple
minimizers of the expected score $Q \mapsto \myS(Q, P)$, the function
$\myS^*$ is well-defined by using a mapping $P \mapsto P^*$ that
chooses a $P^*$ out of the set of minimizers.  If $\myS$ is proper and
$P^* = P$, then $\myS^* = \myS$, so the proper scoring rules are fixed
points under the properization operator.

Importantly, Theorem~\ref{th:properization} is a special case of a
general and powerful construction studied in detail by
\citet{GruenDawid2004} and \citet{Dawid2007}.  Specifically, given
some action space $\cA$ and a loss function $\myL : \cA \times \Omega
\to \Real$, suppose that for each $P \in \cP$ there is a Bayes act
$a_P \in \cA$, such that
\[
\int \myL(a_P,\omega) \dd P(\omega) \leq \int \myL(a, \omega) \dd P(\omega) 
\quad \mbox{for all} \quad a \in \cA. 
\]
Then the function
\[
\myS^* : \cP \times \Omega \to \Real, \quad 
(P,\omega) \mapsto \myS^*(P,\omega) = \myL(a_P,\omega),  
\]
is a proper scoring rule.  Note the natural connection to decision-
and utility-based scoring approaches, where the quality of a forecast
is judged by the monetary utility of the induced acts and decisions
\citep{GranPesa2000, GranMach2006, Ehmetal2016}.

In the remainder of the paper we focus on the above special case in
which the action domain $\cA$ is the class $\cP$.  In
Section~\ref{sec:examples} we identify scattered results in the
literature as prominent special cases of properization
(Examples~\ref{ex:binary}--\ref{ex:PMC}), and we use
Theorem~\ref{th:properization} to construct new proper scoring rules
from improper ones (Examples~\ref{ex:CRPS}--\ref{ex:probScore}).
Section~\ref{sec:existence} gives sufficient conditions for the
existence of Bayes acts and Section~\ref{sec:discussion} contains a
brief discussion. All proofs and technical details are moved to the
\hyperref[sec:appendix]{Appendix}.

\section{Examples}  \label{sec:examples}

This section starts with an example in which we review the ubiquitous
misclassification error from the perspective of properization. We go
on to demonstrate how Theorem~\ref{th:properization} has been used
implicitly to construct proper scoring rules in econometric,
meteorological, and statistical strands of literature.  The notion of
properization simplifies and shortens the respective proofs of
propriety, makes them much more transparent, and puts the scattered
examples into a unifying and principled joint framework.  Further
examples show other facets of properization: The scoring rules
constructed in Example~\ref{ex:CRPS} are original, and the discussion
in Example~\ref{ex:CRPS2} illustrates a connection to the practical
problem of the treatment of observational uncertainty in forecast
evaluation. Finally, Example~\ref{ex:probScore} includes an instance
of a situation in which properization fails.

\begin{example}  \label{ex:binary}
Consider probability forecasts of a binary event, where $\Omega =
\lbrace 0,1 \rbrace$ and $\cP$ is the class of the Bernoulli measures.
We identify any $P \in \cP$ with the probability $p = P(\lbrace 1
\rbrace) \in [0,1]$ and consider the scoring rules
\[
\myS_1(P,\omega) := 1 - p \omega - (1-p)(1-\omega) 
\quad \text{ and } \quad 
\myS_2(P,\omega) := |p-\omega|.
\]
The scoring rule $\myS_1$ corresponds to the mean probability rate
(MPR) in machine learning \citep[p.~30]{Ferrietal2009}.  The scoring
rule $\myS_2$ was first considered by \citet{Dawid1986}. It agrees
with the special case $c_1 = c_2$ in Section~4.2 of
\citet{Parry2016} and corresponds to the mean absolute error (MAE)
as discussed by \citet[p.~30]{Ferrietal2009}.\footnote{As noted by
  \citet{Parry2016}, the improper score $\myS_2$ shares its (concave)
  expected score function $P \mapsto \myS_2(P,P)$ with the proper
  Brier score.  This illustrates the importance of the second
  condition in Theorem 1 of \citet{GneitRaft2007}: For a scoring rule
  $\myS$ the (strict) concavity of the expected score function $ G(P)
  := \myS(P,P)$ is equivalent to the (strict) propriety of $\myS$ only
  if, furthermore, $- \myS(P,\cdot)$ is a subtangent of $- G$ at $P$.}
Both $\myS_1$ and $\myS_2$ are improper with common Bayes act
\[
\textstyle
p^* = \one{p \geq \frac{1}{2}} \in \lbrace 0,1 \rbrace, 
\]
and with the same properized score given by the zero-one rule 
\[
\myS^*(P,\omega) = \left\lbrace 
\begin{array}{cl}
0, & p^* = \omega, \\
1, & \text{otherwise.}
\end{array}
\right.
\]
A case-averaged zero-one score is typically referred to as {\em
  misclassification rate}\/ or {\em misclassification error};
undoubtedly, this is the most popular and most frequently used
performance measure in binary classification.  While the scoring rule
$\myS^*$ is proper it fails to be strictly proper
(\citealp[Example 4]{GneitRaft2007};
\citealp[Section 4.3]{Parry2016}). Consequently, misclassification
error has serious limitations as
a performance measure, as persuasively argued by 
\citet[p.~258]{Harrell2015}, among others.  Nevertheless, the scoring
rule $\myS^*$ is proper, contrary to recent claims of impropriety in
the blogosphere.\footnote{See, e.g., 
  \url{http://www.fharrell.com/post/class-damage/} and
  \url{http://www.fharrell.com/post/classification/}.}
\end{example}

For the remainder of the section, let $\Omega = \real$ and let $\cB$
be the Borel $\sigma$-algebra.  We let $\cL$ be the class of Borel
measures $P$ with a Lebesgue density, $p$.  Furthermore, we write
$\cP_k$ for the measures with finite $k$-th moment and $\cP_k^+$ for
the subclasses when Dirac measures are excluded.  Whenever it
simplifies notation, we identify $P$ with its cumulative distribution
function $x \mapsto P((-\infty, x])$.

\begin{example}  \label{ex:weighted} 
Let $\myS_0$ be a proper scoring rule on some subclass $\cP$ of $\cL$
and let\/ $w$ be a nonnegative weight function such that\/ $0 < \int
w(z) \, p(z) \dd z < \infty$ for $p \in \cP$.  Let
\[ 
\myS : \cP \times \real \to \real, \quad (P,y) \mapsto \myS(P,y) = w(y) \, \myS_0(P,y);    
\] 
this score is improper unless the weight function is constant.
Indeed, by Theorem~1 of \citet{GneitRanjan2011}, the Bayes act $P^*$
under $\myS$ has density
\[
p^*(y) = \frac{w(y) \, p(y)}{\int w(z) \, p(z) \dd z}.
\]
From this we see that the key statement in Theorem~1 of
\citet{HolzKlar2018} constitutes a special case of Theorem
\ref{th:properization}. In the further special case in which $\myS_0$
is the logarithmic score, the properized score~\eqref{eq:S*} recovers
the conditional likelihood score of \cite{Diks2011} up to equivalence,
as noted in Example 1 of \citet{HolzKlar2018}.  For analogous results
for consistent scoring functions see Theorem~5 of \citet{Gneit2011}
and Example~2 of \citet{HolzKlar2018}.
\end{example} 

\begin{example}  \label{ex:error.spread} 
For a probability measure $P \in \cP_4$, let $\mu_P$, $\sigma^2_P$,
and $\gamma_P$ denote its mean, variance, and centered third moment.
Let
\[
\myS(P,y) = \left( \sigma_P^2 - (y - \mu_P)^2 \right)^2 
\]
be the `trial score' in equation~(16) of \citet{Christetal2014}. As
\citet{Christetal2014} show in their Appendix A, any Bayes act
$P^*$ under $\myS$ has mean $ \mu_P + \tfrac{1}{2}
\tfrac{\gamma_P}{\sigma_P^2}$ and variance 
\[
\sigma_P^2 \left( 1 + \frac{1}{4} \frac{\gamma_P^2}{\sigma_P^6} \right), 
\]
so properization yields the spread-error score,
\[
\myS^*(P,y) = \left( \sigma_P^2 - \left( y - \mu_P \right)^2
+ \left( y - \mu_P \right) \frac{\gamma_P}{\sigma_P^2} \right)^2 ,
\]
which is proper relative to the class $\cP_4^+$.  Hence the
construction of the spread-error score in \citet{Christetal2014}
constitutes another special case of Theorem~\ref{th:properization}.
\end{example} 

\begin{example}  \label{ex:PMC} 
The predictive model choice criterion of \citet{LaudIbra1995} and
\citet{GelfGhosh1998} uses the scoring rule $\myS(P,y) = \left( y -
\mu_P \right)^2 + \sigma_P^2$, where $\mu_P$ and $\sigma_P^2$ denote
the mean and the variance of a distribution $P \in \cP_2$,
respectively.  As pointed out by \citet{GneitRaft2007}, this score
fails to be proper.  Specifically, any Bayes act $P^*$ under $\myS$
has mean $ \mu_P$ and vanishing variance, so properization yields the
ubiquitous squared error, $\myS^*(P,y) = \left( y - \mu_P \right)^2$.
\end{example} 

The original scoring rules of Examples~\ref{ex:error.spread}
and~\ref{ex:PMC} can be interpreted as functions $ \myS : \cA \times
\Omega \rightarrow \real$ in the Bayes act setting of
\citet{GruenDawid2004} and \citet{Dawid2007}, where the action space
$\cA$ is given by $\real \times [0, \infty)$.  Hence, the
  properization method can be interpreted as an application of
  Theorem~3 of \citet{Gneit2011} to consistent scoring functions for
  elicitable two-dimensional functionals, as discussed by
  \citet{FissZieg2016}.

Detailed arguments and calculations for the subsequent examples are
deferred to the \hyperref[sec:appendix]{Appendix}.

\begin{example}  \label{ex:CRPS} 
For $\alpha > 0$ consider the scoring rule 
\[
\myS_\alpha (P,y) = \int \left| P(x) - \one{y   \le   x} \right|^\alpha \dd x,
\]
where $P$ is identified with its cumulative distribution function
(CDF).  For $\alpha = 2$ this is the well known proper continuous
ranked probability score (CRPS), as reviewed in Section 4.2 of
\citet{GneitRaft2007}.  For $\alpha = 1$ the score $\myS_\alpha$ was
proposed by \citet{Muelleretal2005}, and \citet{ZamoNaveau2018} show
in their Appendix A that for discrete distributions every Dirac
measure in a median of $P$ is a Bayes act.  The same holds true for
general distributions and for all $\alpha \in (0,1]$.  If $\alpha >
1$, the Bayes act $P^*$ under $\myS_\alpha$ is given by
\begin{equation}  \label{eq:CRPS_alpha_P*}
P^* (x) = \left( 
1 + \left( \frac{1 - P(x)}{P(x)} \right)^{1/(\alpha - 1)} \right)^{-1} \one{P(x) > 0},
\end{equation}
and all in all we see that properization of $\myS_\alpha$ works for
any $\alpha > 0$.

Moreover, in the case $\alpha >1$ the mapping $P \mapsto P^*$ is even
injective. Consequently, if the class $\cP$ is such that $P^* \in \cP$
and $\myS_\alpha(P^*,P)$ is finite for $P \in \cP$, the properized
score~\eqref{eq:S*} is even strictly proper relative to $\cP$.  If
$\alpha \in (1,2]$, this can be ensured by restricting $\myS_\alpha$
to the class $\cP_1$. For $\alpha > 2$ the class $\cP_{\mathrm c}$ of
the Borel measures with compact support is a suitable choice.
\end{example} 

\begin{example}  \label{ex:CRPS2}
\citet[p.~58]{FriedThor2012} propose a modification of the CRPS that
aims to account for observational error in forecast evaluation.
Specifically, they consider the scoring rule
\[
\myS_\Phi (P, y) = \int \left| P(x) - \Phi(x-y) \right|^2 \dd x,
\]
where $\Phi \in \cP_1^+$ represents additive observation error.  This
scoring rule fails to be proper, as for probability measures $P, Q \in
\cP_1$ we have
\begin{equation}  \label{eq:NoisyCRPS}
\myS_\Phi (P, Q) = \crps (P, Q * \Phi) - \crps (\Phi, \Phi),
\end{equation}
where $*$ denotes the convolution operator.  Due to the strict
propriety of the CRPS relative to the class $\cP_1$, the unique Bayes
act under $\myS_\Phi$ is given by $P^* = P * \Phi$.
Theorem~\ref{th:properization} now gives the scoring rule
$\myS(P,y) := \myS_\Phi (P^*, y)$, which is proper relative to
$\cP_1$.

In order to account for noisy observational data in forecast
evaluation, equation~\eqref{eq:NoisyCRPS} suggests using the scoring
rule $\myS (P,y) := \crps (P^*, y)$ if the noise is independent,
additive, and has distribution $\Phi$. This corresponds to predicting
hypothetical true values, to which noise is added before they are
compared to observations. The drawbacks of this approach and
alternative techniques are discussed by~\citet{Ferro2017}. The
associated issues in forecast evaluation remain challenges to the
scientific community at large; see, e.g., \cite{Ebertetal2013}
and \cite{Ferro2017}.

\end{example}

\begin{example}  \label{ex:probScore}
Let $\myS$ be a scoring rule, and let $\Phi \in \cL$ be a distribution
with Lebesgue density $\varphi$.  Suppose $\cP$ is a class of
distributions such that $P * \Phi \in \cP$ for $P \in \cP$.  For $P
\in \cP$ define
\[
\myS^\varphi (P, y) := \int \varphi (x-y) \, \myS (P, x) \dd x,
\]
which is again a scoring rule.  If $\myS$ is proper, a Bayes act
under $\myS^\varphi$ is given by $P^* = P * \Phi$, since $\myS^\varphi
(P, Q) = \myS(P, Q * \Phi )$ for $Q \in \cP$, and if $\myS$ is
strictly proper, the Bayes act is unique. Properization now gives the
proper scoring rule $\myS(P,y) := \myS^\varphi (P^*, y)$. An
interesting special case emerges when substituting the CRPS for
$\myS$. This leads to
\begin{equation}  \label{eq:NoisyCRPS2}
\crps^\varphi (P,y) = \myS_\Phi (P,y ) + \crps (\Phi, \Phi), 
\end{equation}
where $\myS_\Phi$ is the scoring rule in the previous example.  For
another special case, let $c > 0$ and $P \in \cL$, to yield
\[
\mathrm{PS}_c(P, y) := - \int_{y-c}^{y+c} p(x) \dd x, 
\]
which recovers the {\em probability score}\/ of
\citet{Wilsonetal1999}.  We have that $\mathrm{PS}_c = 2c \,
\mathrm{LinS}^{\varphi_c}$, where $ \mathrm{LinS}(P,y) := - p(y)$ is
the improper linear score and $\varphi_c$ is a uniform density on
$[-c,c]$.  Properization is not feasible relative to sufficiently rich
classes $\cP$, as Bayes acts fail to exist under both the linear
score and the probability score.  For details, see the
\hyperref[sec:appendix]{Appendix}.
\end{example}

\section{Existence of Bayes acts}   \label{sec:existence}

In Example~\ref{ex:probScore} we presented a scoring rule that cannot
be properized, due to the non-existence of Bayes acts. This section
addresses the question under which conditions on $\myS$ and $\cP$ a
minimum of the expected score function exists. To illustrate the
ideas, we start with a further example.

\begin{example}  \label{ex:normalized}
Using the notation of Example~\ref{ex:error.spread}, consider the
normalized squared error,
\begin{equation*}
\myS(P,y) = \frac{(y - \mu_P)^2}{\sigma_P^2},  
\end{equation*}
as a scoring rule on the classes $\cP_{2,m}$ of the Borel measures
with variance at most $m$, and $\cP_2 = \cup_{m > 0} \cP_{2,m}$,
respectively.  Relative to $\cP_{2,m}$ any Bayes act $P^*$ under
$\myS$ has mean $\mu_P$ and variance $m$, so properization yields
(non-normalized) squared error up to equivalence.  Relative to $\cP_2$
however, there is no Bayes act, since increasing the variance will
always lead to a smaller expected score.
\end{example}

We now turn to a general perspective and discuss sufficient
conditions for the existence of Bayes acts. At first, consider a
finite probability space $\Omega = \lbrace \omega_1, \ldots, \omega_k
\rbrace$. In this situation, geometrical arguments yield  sufficient
conditions. In particular, a Bayes act under $\myS$ exists if the
\textit{risk set}
\[
\mathcal{S} := \lbrace (x_1, \ldots, x_k) \mid \exists \, P \in
\mathcal{P} : x_j = \myS (P, \omega_j) , \, j=1,\ldots, k \rbrace
\subset \mathbb{R}^k
\]
is closed from below and bounded from below; see Theorem~1 in
Chapter~2.5 of \citet{Ferguson1967}. Extending this result to a general
sample space $\Omega$ is non-trivial since in this case $\mathcal{S}$
can be a subset of an infinite-dimensional vector space. In the
following we employ well-known concepts of functional analysis in order
to discuss two possible extensions. All proofs are deferred to the
\hyperref[sec:appendix]{Appendix}.

Let $\cP$ be a set of probability measures on a general probability
space $\Omega$ and let $\cA$ be a topological space. We return to the
setting of Section~\ref{sec:introduction} and consider functions
$\myS: \cA \times \Omega \rightarrow \real$. This makes the results
more general and easier to apply in situations where the scoring rule
depends on $P$ only via some finite number of parameters. Concerning
the latter point, note that the normalized squared error of
Example~\ref{ex:normalized} can be written as a composition of the
mapping $P \mapsto (\mu_P, \sigma^2_P)$ and the function $s(x_1,x_2,y)
:= (y -x_1)^2 / x_2$, with $s$ being defined on $\real \times (0,
\infty) \times \real$. Consequently, the expected normalized squared
error attains its minimum if the expected score of $s$ attains its
minimum. Note that such a decomposition of the scoring rule is
possible for Examples~\ref{ex:error.spread} and~\ref{ex:PMC} as well,
as alluded to in the comments that succeed these examples.  

We impose the following integrability assumption on $\myS$.

\begin{definition}
The mapping $\myS: \cA \times \Omega \rightarrow \real$ is
\textit{uniformly bounded from below} if there exists a function
$g: \Omega \rightarrow \real$ which is integrable with respect to any
$P \in \cP$ and such that $\myS (a,\cdot) \geq g$ holds for all
$a \in \cA$.
\end{definition}

Our first result is similar to Theorem~2 in Chapter~2.9 of
\citet{Ferguson1967}, which proves the existence of minimax decision
rules. 

\begin{theorem}  \label{th:existence1}
Suppose $\myS$ is lower semicontinuous in its first component and
uniformly bounded from below. If\/ $\cA$ is compact, then the function
$a \mapsto \myS(a,P)$ attains its minimum for any $P \in \cP$.
\end{theorem}

This theorem can be used to prove the existence of a Bayes act
for a given scoring rule. However, it is not applicable to
Example~\ref{ex:normalized}. To see this, recall the decomposition of
$\myS$ mentioned above and note that restricting $\myS$ to $\cP_{2,m}$
corresponds to restricting $s$ to $\real \times (0,m]$. The latter set
is not a compact space and neither is its closure. Consequently, we aim
to dispense with the compactness assumption used in
Theorem~\ref{th:existence1}.

To do so, we need additional concepts from functional analysis. Let
$\cX$ be a real normed vector space.  Recall that a function $h : \cX
\rightarrow \real$ is called \textit{coercive} if for any sequence
$(x_n)_{n \in \mathbb{N}} \subset \cX$ the implication
\begin{equation*}
\lim_{n \rightarrow \infty} \Vert x_n \Vert = \infty  \quad \Rightarrow \quad \lim_{n \rightarrow \infty} h(x_n) = \infty
\end{equation*}
holds true, see, e.g.~Definition~III.5.7 in~\citet{Werner2018}. By
\textit{weak topology} on $\cX$, we mean the weakest topology such
that all real-valued linear mappings on $\cX$ are continuous; see,
e.g.~Chapters~2.13 and~6.5 in~\citet{AlipBord2006}. The space $\cX$ is
called a \textit{reflexive Banach space} if it is complete and the
canonical embedding of $\cX$ into its bidual space is surjective; see,
e.g.~Chapter~III.3 in~\citet{Werner2018} or Chapter~6.3
in~\citet{AlipBord2006}. Combining these concepts, we obtain a
complement to Theorem~\ref{th:existence1}.

\begin{theorem}  \label{th:existence2}
Let $\cA$ be a weakly closed subset of a reflexive Banach space.
Moreover, suppose $\myS$ is weakly lower semicontinuous in its first
component and uniformly bounded from below. If the function $a \mapsto
\myS(a,P)$ is coercive, then it attains its minimum.
\end{theorem}

This result yields the existence of Bayes acts as long as the
integrated scoring rule is coercive for any $P \in \cP$, where $\cP$
is a reflexive Banach space. To conclude this section, we connect
Theorem~\ref{th:existence2} to Example~\ref{ex:normalized}: The
function $s(\cdot, \cdot, y)$ from the decomposition of $\myS$
mentioned above is defined on $\real \times (0, \infty)$, which is a
subset of the reflexive Banach space $\real^2$. Moreover, $s$ is
bounded from below by zero and continuous in its first component. As
mentioned above, restricting the class $\cP_2$ to $\cP_{2,m}$
corresponds to restricting the domain of $s$ to $\real \times (0, m]$
and in this situation, integrating $s$ with respect to $y$ gives a
coercive function. Consequently, Theorem~\ref{th:existence2} can be
used to show that $\myS$ can be properized if restricted to
$\cP_{2,m}$.

\section{Discussion}   \label{sec:discussion}

In this article we have introduced the concept of properization, which
is rooted in the Bayes act construction of \citet{GruenDawid2004} and
\citet{Dawid2007}, and we have drawn attention to its widespread
implicit use in the transdisciplinary literature on proper scoring
rules, where our unified approach yields simplified, shorter, and
considerably more instructive and transparent proofs than extant
methods.  Moreover, using new examples, we have demonstrated the power
of the properization approach in the creation of new proper scoring
rules from existing improper ones.  

Since the central element in the construction of a properized score is
a Bayes act, we have discussed conditions on the scoring rule $\myS$
and the class $\cP$ that guarantee its existence.  Undoubtedly, there
are alternative paths to existence results in the spirit of
Theorems~\ref{th:existence1} and~\ref{th:existence2}, and the
derivation of sufficient conditions in alternative situations is an
interesting open problem.  Furthermore, we have not explored necessary
conditions for the existence of Bayes acts in this work. Their
derivation and the refinement of sufficient conditions on $\myS$ and
$\cP$ remain challenges that we leave for future work.

\section*{Appendix: Proofs} \label{sec:appendix}
\addcontentsline{toc}{section}{Appendix: Proofs}

Here we present detailed arguments for the technical claims in
Examples~\ref{ex:CRPS}, \ref{ex:CRPS2}, and~\ref{ex:probScore}
as well as the proofs of Theorems~\ref{th:existence1}
and~\ref{th:existence2}.

\subsection*{Details for Example~\ref{ex:CRPS}}
\addcontentsline{toc}{subsection}{Details for Example~\ref{ex:CRPS}}

We fix some distribution $P$ and start with the case $\alpha > 1$. An
application of Fubini's theorem gives
\begin{equation} \label{eq:Ex4SQP}
\myS_\alpha (Q,P) = \int \int \vert Q(x) - \one{y \leq x} \vert^\alpha  \dd P(y) \dd x .
\end{equation}
Given $x \in \real$ we seek the value $Q(x) \in [0,1]$ that minimizes
the inner integral in~\eqref{eq:Ex4SQP}.  If $x$ is such that $P(x)
\in \lbrace 0, 1 \rbrace$, the equality $\one{y \leq x} = P(x)$ holds
for $P$-almost all $y$, hence $Q(x)= P(x)$ is the unique minimizer.
If $x$ satisfies $P(x) \in (0,1)$, define the function
\[
g_{x,P}(q) := \int \vert q - \one{y \leq x} \vert^\alpha \dd P(y) 
= (1-P(x)) q^\alpha + P(x) (1-q)^\alpha,
\]
which is strictly convex in $q \in (0,1)$ with derivative
\[
g_{x,P}'(q) = \alpha (1- P(x)) q^{\alpha - 1} - \alpha P(x) (1- q)^{\alpha - 1}
\]
and a unique minimum at $q = q_{x,P}^* \in (0,1)$.  As a consequence,
the minimizing value $Q(x)$ is given by
\begin{equation*}
Q(x) = q_{x,P}^* = \left( 1 + \left( \frac{1-P(x)}{P(x)} \right)^{1/(\alpha - 1)} \right)^{-1}.
\end{equation*}
The function $Q$ defined by the minimizers $Q(x)$, $x \in \real$ is a
minimizer of $\myS_\alpha ( \cdot ,P)$ and if $\myS_\alpha (Q, P)$ is
finite, it is unique Lebesgue almost surely.  Since $\alpha >1$ the
function $Q$ has the properties of a distribution function and hence
$P^*$ defined by~\eqref{eq:CRPS_alpha_P*} is a Bayes act for
$P$. Moreover, equation~\eqref{eq:CRPS_alpha_P*} shows that the
relation between $P$ and $P^*$ is one-to-one.

It remains to be checked under which conditions the properization of
$\myS_\alpha$ is not only proper but strictly proper.  The
representation~\eqref{eq:CRPS_alpha_P*} along with two Taylor
expansions implies that $P^*$ behaves like $P^{1/(\alpha -1)}$ in the
tails. This has two consequences. At first, the above arguments show
that for $\myS_\alpha (P^*, P)$ to be finite $x \mapsto g_{x,P}
(P^*(x))$ has to be integrable with respect to Lebesgue measure.
Hence, the tail behavior of $P^*$ and the inequality $\alpha/(\alpha
- 1) > 1$ for $\alpha > 1$ show that $\myS_\alpha (P^*, P)$ is finite
for $P \in \cP_1$.  Second, $P^*$ has a lighter tail than $P$ for
$\alpha \in (1,2)$ and a heavier tail for $\alpha > 2$.  In the latter
case $P \in \cP_1$ does not necessarily imply $P^* \in \cP_1$.  Hence,
without additional assumptions, strict propriety of the properized
score~\eqref{eq:S*} can only be ensured relative to $\cP_\mathrm{c}$
for $\alpha > 2$ and relative to the class $\cP_1$ for $\alpha \in (1,
2]$.
 
We now turn to $\alpha \in (0,1)$.  In this case, the function
$g_{x,P}$ is strictly concave, and its unique minimum is at $q = 0$
for $P(x) < \frac{1}{2}$ and at $q = 1$ for $P(x) > \frac{1}{2}$.  If
$P(x) = \frac{1}{2}$, then both $0$ and $1$ are minima.  Arguing as
above, every Bayes act $P^*$ is a Dirac measure in a median of $P$.

Finally, $\alpha = 1$ implies that $g_{x,P}$ is linear, thus, as for
$\alpha \in (0,1)$, every Dirac measure in a median of $P$ is a Bayes
rule.  The only difference to the case $\alpha \in (0,1)$ is that if
there is more than one median, there are Bayes acts other than Dirac
measures, since $g_{x,P}$ is constant for all $x$ satisfying $P(x) =
\frac{1}{2}$.

\subsection*{Details for Example~\ref{ex:CRPS2}}
\addcontentsline{toc}{subsection}{Details for Example~\ref{ex:CRPS2}}

Let $P, Q$ and $\Phi$ be distribution functions.  By the
definition of the convolution operator
\begin{equation*}
\int \one{y \leq x}  \dd (Q * \Phi) (y) = \int \Phi (x-y)  \dd Q(y) 
\end{equation*}
holds for $x \in \real$.  Using this identity and Fubini's theorem
leads to
\begin{align*}
\myS_\Phi (P,Q) 
& = \int \! \int \left( P(x)^2 - 2 P(x) \Phi(x-y) + \Phi(x-y)^2 \right) \dd Q(y) \dd x \\
& = \int \! \int \left( P(x)^2 - 2 P(x) \one{y \leq x} + \one{y \leq x} \right) \dd (Q * \Phi)(y) \dd x \\
& \quad + \int \! \int \Phi(x-y) (\Phi(x-y) - 1) \dd Q(y) \dd x \\
& = \int \! \int (P(x) - \one{y \leq x})^2  \dd x \dd (Q * \Phi)(y) 
  - \int \Phi (x) (1- \Phi (x))  \dd x,
\end{align*}
which verifies equality in~\eqref{eq:NoisyCRPS}.  Moreover, the strict
propriety of the CRPS relative to the class $\cP_1$ gives $\myS_\Phi
(P, Q) < \infty$ for $P, Q, \Phi \in \cP_1$, thereby demonstrating
that the Bayes act is unique in this situation.

\subsection*{Details for Example~\ref{ex:probScore}}
\addcontentsline{toc}{subsection}{Details for Example~\ref{ex:probScore}}

For distributions $P, Q \in \cP$ and $c > 0$, the Fubini-Tonelli
theorem and the definition of the convolution operator give
\begin{align*}
\myS^\varphi (P,Q) &= - \int \int \varphi (x-y) \myS(P,x) \dd Q(y) \dd x \\
&= \int \int \varphi (x-y) \dd Q(y) \, \myS(P,x) \dd x = \myS (P, Q * \Phi),
\end{align*}
so the stated (unique) Bayes act under $\myS^\varphi$ follows from
the (strict) propriety of $\myS$.  Proceeding as in the details for
Example~\ref{ex:CRPS2} we verify identity~\eqref{eq:NoisyCRPS2}.

For $P \in \cL$ the same calculations as above show that the
probability score satisfies
\[
\mathrm{PS}_c(P,Q) = 
2c \int \frac{Q(x + c) - Q(x - c)}{2c} \: \mathrm{LinS}(P,x) \dd x,
\]
where $\mathrm{LinS}(P,y) = - p(y)$ is the linear score.
Consequently, to demonstrate that Theorem~\ref{th:properization} is
neither applicable to $\mathrm{PS}_c$ nor to $\mathrm{LinS}$, it
suffices to show that there is a distribution $Q$ such that $P \mapsto
\mathrm{LinS}(P,Q)$ does not have a minimizer.  We use an argument
that generalizes the construction in Section 4.1 of
\citet{GneitRaft2007} who show that $\mathrm{LinS}$ is improper.  Let
$q$ be a density, symmetric around zero and strictly increasing on
$(-\infty, 0)$.  Let $\epsilon > 0$ and define the interval $I_k :=
((2k - 1) \epsilon, (2k + 1) \epsilon]$ for $k \in \mathbb{Z}$.
  Suppose $p$ is a density with positive mass on some interval $I_k$
  for $k \neq 0$.  Due to the properties of $q$, the score
  $\mathrm{LinS}(P,Q)$ can be reduced by substituting the density
defined by 
\[
\tilde{p}(x) := p(x) - \one{x \in I_k} \, p(x) + \one{x + 2k \epsilon
  \in I_k} \, p(x + 2k \epsilon)
\]
for $p$, i.e., by shifting all probability mass from $I_k$ to the
modal interval $I_0$.  Repeating this argument for any $\epsilon > 0$
shows that no density $p$ can be a minimizer of the expected score
$\mathrm{LinS}(P,Q)$.  Note that the assumptions on $q$ are stronger
than necessary in order to facilitate the argument.  They can be
relaxed at the cost of a more elaborate proof.

\subsection*{Proof of Theorem~\ref{th:existence1}}
\addcontentsline{toc}{subsection}{Proof of Theorem~\ref{th:existence1}}

Let $(a_n)_{n \in \mathbb{N}} \subset \cA$ be a sequence with $ a :=
\lim_{n \rightarrow \infty} a_n$. Since $\myS$ is lower semicontinuous
in its first component and uniformly bounded from below by $g$, Fatou's
lemma gives
\begin{equation*}
\liminf_{n \rightarrow \infty} \int \myS (a_n, \omega) \dd P(\omega)
\geq \int \liminf_{n \rightarrow \infty} \myS (a_n,\omega) \dd
P(\omega) \geq \myS(a,P)
\end{equation*} 
for any $P \in \cP$. Hence, $a \mapsto \myS(a, P)$ is a lower
semicontinuous function for any $P \in \cP$ and due to the assumed
compactness of $\cA$, the result now follows from Theorem~2.43
in~\citet{AlipBord2006}.

\subsection*{Proof of Theorem~\ref{th:existence2}}
\addcontentsline{toc}{subsection}{Proof of Theorem~\ref{th:existence2}}

The same arguments as in the proof of Theorem~\ref{th:existence1}
show that $a \mapsto \myS(a, P)$ is a weakly lower semicontinuous
function for any $P \in \cP$. If $P \in \cP$ is such that this function
is also coercive, then proceeding as in the proof of Satz~III.5.8
in~\citet{Werner2018} gives a weakly convergent sequence
$(a_n)_{n \in \mathbb{N}} \subset \cA$ with $\lim_{n \rightarrow
\infty} \myS( a_n, P) = \inf_{a \in \cA} \myS(a, P)$. Since $\cA$ is
weakly closed by assumption, it contains the weak limit  $a^*$ of the
sequence $(a_n)_{n \in \mathbb{N}}$ and hence weak lower semicontinuity
implies that $a \mapsto \myS(a, P)$ attains its minimum at $a^* \in
\cA$.

\section*{Acknowledgments}
\addcontentsline{toc}{section}{Acknowledgments}

Tilmann Gneiting is grateful for funding by the Klaus Tschira
Foundation and by the European Union Seventh Framework Programme under
grant agreement 290976.  Part of his research leading to these results
has been done within subproject C7 ``Statistical postprocessing and
stochastic physics for ensemble predictions'' of the Transregional
Collaborative Research Center SFB / TRR 165 ``Waves to Weather''
(\url{www.wavestoweather.de}) funded by the German Research Foundation
(DFG).  Jonas Brehmer gratefully acknowledges support by DFG through
Research Training Group RTG 1953.  We thank Tobias Fissler and Matthew
Parry for instructive discussions.

\end{document}